\documentclass[11pt]{amsart}
\usepackage{amssymb}
\usepackage{mathrsfs}
\usepackage{hyperref}
\usepackage{color}
\usepackage{graphicx}
\usepackage{epstopdf}
\usepackage{fancybox}
\usepackage{pstricks}
\usepackage{url}
\usepackage{mathtools}
\usepackage{tensor}

\date{\today}

\newcommand{\nn}{\nonumber}
\newcommand{\be}{\begin{equation}}
\newcommand{\ee}{\end{equation}}

\DeclareMathOperator{\loc}{loc}

\newcommand{\om}{\omega}

\def\lb{\label}
\def\qq{\qquad}

\def\Proof{\noindent{\bf Proof} \quad}
\def\qed{\hfill $\Box$ \smallskip}

\def\om{\omega}
\def\loc{{\rm loc}}
\def\bb{\begin}
\def\la{\langle}
\def\ra{\rangle}

\def\bea{\begin{eqnarray}}
\def\eea{\end{eqnarray}}
\def\beaa{\begin{eqnarray*}}
\def\eeaa{\end{eqnarray*}}

\newtheorem{thm}{Theorem}[section]
\newtheorem{rem}[thm]{Remark}
\newtheorem{lem}[thm]{Lemma}

\numberwithin{equation}{section}

\begin{document}

\title[]{Asymptotic stability of the spectra of generalized indefinite strings}

\author[G. Xu]{Guixin Xu$^*$}

\address{School of Mathematics and Statistics, Beijing Technology and Business University, Beijing, China}
\email{\href{mailto:guixinxu$_-$ds@163.com}{guixinxu$_-$ds@163.com}}

\author[M. Zhang]{Meirong Zhang}

\address{Department of Mathematical Sciences, Tsinghua University,
Beijing 100084, China}
\email{\href{mailto:zhangmr@tsinghua.edu.cn}{zhangmr@tsinghua.edu.cn}}

\author[Z. Zhou]{Zhe Zhou}

\address{Academy of Mathematics and Systems Science,  Chinese Academy of Sciences,  Beijing 100190, China}
\address{School of Mathematical Sciences, University of Chinese Academy of Sciences, Beijing 100049, China}
\email{\href{mailto:zzhou@amss.ac.cn}{zzhou@amss.ac.cn}}

\thanks{$^*$ Corresponding author}
\thanks{M.\ Z.\ was supported by the National Natural Science Foundation of China (Grant No. 11790273)}
\thanks{Z.\ Z.\ was supported by the National Natural Science Foundation of China (Grant No. 12271509,12090010 and 12090014)}



\begin{abstract}
This paper focuses on the asymptotic stability of the spectra of generalized indefinite strings (GISs). A unitarily equivalent linear relation is introduced for GISs. It is shown that the solutions of the corresponding differential equations are continuously dependent on distributions and measures under certain conditions. Using these results, the convergence of unitarily equivalent linear relations for GISs is discussed. By the perturbation theory to closed linear relations, an asymptotic stability result concerning the spectra of linear relations for GISs is obtained.
\end{abstract}

\maketitle


\section{Introduction}

A generalized indefinite string (GIS, for short) is a triple $(L,\om,\nu)$ such that $L\in(0,\infty]$, $\om$ is a real-valued distribution in $H_\loc^{-1}[0,L)$ and $\nu$ is a non-negative Borel measure on $[0,L)$.  Associated with such a triple is the ordinary differential equation of the form
\be\lb{gis}
-f''=z\om f+z^2\nu f.
\ee
on $[0,L)$, where $z$ is a complex spectral parameter. The spectral problems of this equation have attracted significant attention from both mathematicians and physicists. Firstly, they constitute a canonical model for operators with a simple spectrum \cite{EckhardtK16, EckhardtK18}. Secondly, they play an important role in the connection with certain completely integrable nonlinear wave equations, such as the Camassa-Holm equation \cite{CamassaH93}, where these spectral problems arise as isospectral problems \cite{Dou23, EckhardtK14, EckhardtT13}.

J. Eckhardt and his collaborators have done extensive work on GISs \cite{EckhardtK16, EckhardtK18, EckhardtK21, EckhardtK24, EckhardtKN19}. They established in \cite{EckhardtK16} that there exists a homeomorphism between the set of all GISs and the set of all Herglotz-Nevanlinna functions. Very recently, they obtained a couple of Szeg\H{o}-type theorems for this correspondence in \cite{EckhardtK24}. Utilizing these results, the stability results for the absolutely continuous spectra of some model examples of GISs are established \cite{EckhardtK21, EckhardtKK22}. It is noteworthy that the spectrum of a GIS is defined as the spectrum of its corresponding linear relation and some stability results concerning the spectra of linear relations under certain perturbations have been obtained in \cite{XS17, XS18}. The aim of this paper is to discuss the further problem: {\it What is the behaviour of the spectra of $T$  if the distribution $\om$ and Borel measure $\nu$ change}?

It is worth mentioning that G. Meng and M. Zhang \cite{Meng13} studied the dependence of solutions and eigenvalues of second-order linear measure differential equations on measures, obtaining two fundamental results regarding continuity and continuous Fr\'echet differentiability in measures under the weak$^*$ topology and the norm topology of total variations, respectively. Recently, J. Chu and his collaborators \cite{Chu20} proved that the eigenvalues for periodic generalized Camassa-Holm equations are continuous in potentials under the weak topology and provided some estimates for them.

In the present paper, we consider the regular case of GISs \cite{EckhardtK18}. Without loss of generality, we assume that $L=1$. Note that the associated relation $T$ of a GIS $(1,\om,\nu)$ belongs to the Hilbert space ${\mathcal H}^2$. For the definition of $\mathcal{H}$, see \eqref{space}. Since ${\mathcal H}$ depends on $\nu$, we can not directly discuss the relationship between the perturbation of $T$ and the variations of $\om$ and $\nu$. The key idea here is to introduce a new Hilbert space $\widetilde{\mathcal{H}}$ and a new linear relation $\widetilde{T}$, which is unitarily equivalent to $T$ under some conditions. We then study the relationship between the perturbation of $\widetilde{T}$ and the variations of $\om$ and $\nu$. Together with the stability results concerning the spectra of linear relations in \cite{XS17, XS18}, this implies the perturbation results of the spectrum of $\widetilde{T}$. Therefore, the perturbation results of the spectrum of $T$ can be obtained since unitarily equivalent relations have the same spectral structure.

The rest of this paper is organized as follows. In Section 2, we recall some basic concepts and fundamental results on the self-adjoint relation $T$ associated with a GIS, and introduce a unitary operator and a unitarily equivalent relation $\widetilde{T}$ of $T$. Additionally, the relationship between the generalized convergence and the norm resolvent convergence of relations is established. In Section 3, we provide properties of solutions of the corresponding differential equation and discuss the continuous dependence of solutions. For $n\geq1$, let $T_n$ be the self-adjoint relations of  $(1,\om_n,\nu_n)$ and $\widetilde{T}_n$ be their unitarily equivalent relations. We prove that $\{\widetilde{T}_n\}$  converges to $\widetilde{T}$ in the sense of norm resolvent convergence when $\{\om_n\}$ and $\{\nu_n\}$ converge to $\om$ and $\nu$ in some sense, respectively. By applying the small gap perturbation result for the spectra of closed linear relations in \cite{XS17}, we obtain an asymptotic stability result about the spectra of linear relations for GISs on distributions and measures.

\medskip

%
%
%
%
%

\section{Preliminaries}\label{sec.2}

In this section, we recall some notations and basic concepts, and give some fundamental results about generalized indefinite strings (GISs). In particular, we give some properties of Borel measures and self-adjoint relations associated with GISs.

%
%
%
%
%

\subsection{Generalized indefinite strings}

Let $L=1$, we consider the regular case of GISs. First we introduce the following spaces
$$\begin{array}{rrll}
	H_\loc^1[0,1)&:=&\{f\in AC_\loc[0,1):\, f'\in L_\loc^2[0,1)\},\\[0.8ex]
	H^1[0,1)&:=&\{f\in H_\loc^1[0,1):\, f,\,f'\in L^2[0,1)\},\\[0.8ex]
	H_{c}^1[0,1)&:=&\{f\in H^1[0,1):\, \rm{supp}(f) \mbox{ is compact in }\,[0,1)\}.\\[0.6ex]
\end{array}$$
The space of distributions $H_\loc^{-1}[0,1)$ is the topological dual of $H_{c}^1[0,1)$. In fact, any $\mathbf{q}\in L_\loc^2[0,1)$ corresponds to a distribution
$\chi\in H_\loc^{-1}[0,1)$,
defined by
\be
\chi(h):=-\int_{0}^{1}\mathbf{q}(x)h'(x)dx,\qq h\in H_{c}^1[0,1).\nn
\ee
Moreover, the mapping $\mathbf{q}\rightarrow \chi$
is a one-to-one correspondence between $L_\loc^2[0,1)$ and $H_\loc^{-1}[0,1)$.
The function $\mathbf{q}\in L_\loc^2[0,1)$ is referred to as the {\it normalized anti-derivative} of the distribution $\chi\in H_\loc^{-1}[0,1)$.
We say that $\chi$ is real-valued if $\mathbf{q}$ is real-valued almost everywhere on $[0,1)$.

A special class of distributions in $H_\loc^{-1}[0,1)$ consists of all Borel measures on $[0,1)$. When $\chi$ is a complex-valued Borel measure on $[0,1)$, the corresponding distribution in $H_\loc^{-1}[0,1)$ is given by the following linear functional
\be
h \mapsto \int_{[0,1)}hd\chi,\qq h\in H_{c}^1[0,1),\nn
\ee
whose normalized anti-derivative $\mathbf{q}$ is the following left-continuous function
\be
\mathbf{q}(x)=\int_{[0,x)}d\chi, \qq \mbox{a.e. } x\in [0,1).
\ee

Associated with a GIS $(1,\om,\nu)$ is the differential equation
\be\lb{GIS}
-f''=z\om f+z^2\nu f,
\ee
where $z$ is a complex spectral parameter. This differential equation has to be understood in a weak sense. A {\it solution} of \eqref{GIS} is a function $f\in H_\loc^1[0,1)$ such that
\be
\Delta_fh(0)+\int_0^1f'(x)h'(x)dx=z\om(fh)+z^2\int_{[0,1)}f(t)h(t)d\nu(t), \nn
\ee
for all $h\in H_c^1[0,1)$ and some constant $\Delta_f\in {\mathbb C}$. In this case, the constant $\Delta_f$ is denoted by $f'(0-)$ via the integration by parts formula.

Let $\mathbf{u}$  be the normalized anti-derivatives of the distributions $z\om +z^2\nu$. Then a function $f\in H_\loc^1[0,1)$ is a solution of \eqref{GIS} if and only if
\be\lb{solution2}
f'(x)+\mathbf{u}(x)f(x)=\Delta_f+\int_0^x\mathbf{u}(t)f'(t)dt,
\ee
for almost everywhere $x\in[0,1)$. Here $\Delta_f\in {\mathbb C}$ is the same as in the definition above.

The Wronski determinant $W(\theta, \phi)$ of two solutions $\theta, \phi$ of \eqref{GIS} is defined by
\be
W(\theta, \phi)=\theta(0)\phi'(0-)-\theta'(0-)\phi(0).\nn
\ee
For more details on the Wronski determinant, see \cite[Section 3]{EckhardtK16}.

From now on, for any GIS $(1,\om,\nu)$, we always assume that $\nu$ satisfies the following basic hypothesis unless otherwise stated.

\medskip
\noindent{\bf Hypothesis (0).}
{\it $\nu$ is a positive and finite Borel measure on $[0,1)$, and absolutely continuous with respect to the Lebesgue measure.}
\medskip

By \cite[Theorem 3.3.13]{Yan04}, there exists an almost everywhere positive measurable function $p(x)$ such that
\be\lb{ec1}
\int_{[0,1)}fdv=\int_{0}^1(fp)dx,\;f\in L^1([0,1);\nu).
\ee
The function $p$ is called the {\it Radon-Nikodym derivative} of $\nu$ with respect to the Lebesgue measure. Define an operator $P: L^2([0,1);\nu)\to  L^2[0,1)$ by
\be
Pf:=f\sqrt{p},\;\; f\in L^2([0,1);\nu),
\ee
where $p$ is the Radon-Nikodym derivative of $\nu$. Obviously, the operator $P$ is well-defined by \eqref{ec1}.

\bb{lem}\lb{lc1}
$P$ is a unitary operator from $L^2([0,1);\nu)$ to $L^2[0,1)$.
\end{lem}

\Proof
For any $f\in L^2([0,1);\nu)$, by \eqref{ec1} we have
\be
\|Pf\|_{L^2[0,1)}^2=\int_0^1|f\sqrt{p}|^2dx=\int_0^1|f|^2pdx=\int_{[0,1)}|f|^2dv=\|f\|_{L^2([0,1);\nu)}^2,\nn
\ee
where $p$ is the Radon-Nikodym derivative of $\nu$. This implies that $P$ is an isometry.

For any $h\in L^2[0,1)$, let  $f:=\frac{h}{\sqrt{p}}$. Again by \eqref{ec1}, we have
\be
\int_{[0,1)}|f|^2dv=\int_{[0,1)}\frac{|h|^2}{p}dv=\int_0^1|h|^2dx<\infty.\nn
\ee
This implies that $f\in L^2([0,1);\nu)$, and $Pf=h$. Then $P$ is surjective. Thus, $P$ is a unitary operator from $L^2([0,1);\nu)$ to $L^2[0,1)$. The proof is complete.
\qed

For $n\geq1$, we assume that $p$ and $p_n$ are the Radon-Nikodym derivatives of $\nu$ and $\nu_n$. The following hypothesis on Radon-Nikodym derivatives is  central to this paper.

\medskip
\noindent{\bf Hypothesis (1).}
{\it There exists a non-negative integrable function g on $[0,1)$ such that $p_n\leq g$, a.e., and $p_n\xrightarrow[]{a.e.}p$ as $n\to\infty$.}
\medskip

The following result can be directly derived from Lebesgue's Dominated Convergence Theorem and \cite[Theorem 3.4.6]{Yan04}. We omit its proof.

\bb{lem}\lb{lc3}
Assume that  Hypothesis {\rm(1)} holds. Then
\be
\int_0^1\left|\sqrt{p_n}-\sqrt{p}\right|^2dx\to0,\;\int_0^1\left|p_n-p\right|dx\to0\;\;as\;\;n\to\infty.\nn
\ee
\end{lem}

The following hypothesis on distributions is also central to this paper. First we introduce a definition of the convergence of distributions.

\bb{defn} \label{dis-uc}
{\rm For $n\geq1$, let $\om_{n},\ \om\in H_\loc^{-1}[0,L)$. We say that $\om_{n}$ converges to $\om$ in the sense of anti-derivative $L^2$ convergence if ${\mathbf w_n}$ converges to ${\mathbf w}$ in $L^2[0,1)$, that is,
	\be\lb{e5}
	\int_0^1|{\mathbf w_n}-{\mathbf w}|^2dx\to0,\;\; \mbox{as}\;\;n\to\infty,
	\ee
	where ${\mathbf w}$ and ${\mathbf w_n}$ are the anti-derivatives of $\om$ and $\om_n$.}
\end{defn}

\medskip
\noindent{\bf Hypothesis (2).}
{\it $\om_{n}$ converges to $\om$ in the sense of anti-derivative $L^2$ convergence.}
\medskip

%
%
%
%
%

\subsection{Self-adjoint linear relations associated with generalized indefinite strings}

To study the GIS $(1,\om,\nu)$, we recall from \cite{EckhardtK16} the following space
\be
H'[0,1):=\{f\in H_\loc^1[0,1):\,f'\in L^2[0,1),\;\displaystyle\lim_{x\to 1^{-}}f(x)=0\}  \nn
\ee
as well as the linear space
\be\lb{space}
\mathcal{H}:=H'_{*}[0,1)\times L^2([0,1);\nu).
\ee
Here
\be
H'_{*}[0,1)=\{f\in H'[0,1):\,f(0)=0\}.\nn
\ee
Note that $\mathcal{H}$ is a Hilbert space equipped with the scalar product
\be
\la f,g\ra_{\mathcal{H}}:=\int_0^1f'_1(x)g'_1(x)^*dx+\int_{[0,1)}f_2(x)g_2(x)^*d\nu(x).\nn
\ee

For any pair $(f,g)\in \mathcal{H} \times \mathcal{H}$, where $f=(f_1,f_2)\in\mathcal{H}$ and $g=(g_1,g_2)\in \mathcal{H}$, we use the following equation
    \be\lb{relation}
-f''_1=\om g_1+\nu g_2,\qq \nu f_2=\nu g_1,
    \ee
to define a linear subspace $T$ of $\mathcal{H} \times \mathcal{H}$. It is also called a linear relation or a multi-valued linear operator on $\mathcal{H}$ and is denoted by $T \in LR(\mathcal{H})$. Moreover, $T$ is a self-adjoint linear relation in $\mathcal{H}$; see \cite[Theorem 4.1]{EckhardtK16}.

Let
\be\lb{space1}
\widetilde{\mathcal{H}}:=H'_{*}[0,1)\times L^2[0,1).
\ee
Then $\widetilde{\mathcal{H}}$ is a Hilbert space equipped with the inner product
\be
\la f,g\ra_{\widetilde{\mathcal{H}}}:=\int_0^1f'_1(x)g'_1(x)^*dx+\int_0^1f_2(x)g_2(x)^*dx,\nn
\ee
where $f=(f_1,f_2)\in\widetilde{\mathcal{H}}$ and $g=(g_1,g_2)\in \widetilde{\mathcal{H}}$. Define the operator $Q:\mathcal{H}\to \widetilde{\mathcal{H}}$ by
\be
Qf:=(f_1,Pf_2),\;\;f=(f_1,f_2)\in\mathcal{H}.
\ee
By Lemma \ref{lc1}, we have the following result.

\bb{lem}\lb{lc2}
$Q$ is a unitary operator from $\mathcal{H}$ to $\widetilde{\mathcal{H}}$.
\end{lem}

Let $\widetilde{T}:=QTQ^{-1}=QTQ^*$. Then $\widetilde{T}$ is a self-adjoint linear relation in $\widetilde{\mathcal{H}}$. In particular, $\widetilde{T}$ and $T$ have the same spectral structure. By \cite[Proposition 4.3]{EckhardtK16}, we obtain the following result.

\bb{lem}\lb{l3}
If $z$ belongs to the resolvent set $\rho(\widetilde{T})$, then
\be
z(\widetilde{T}-z)^{-1}g(x)=Q\left(\langle g,Q\mathcal{G}(x,\cdot)^*\rangle_{\widetilde{\mathcal H}}\left(
                                                                      \begin{array}{c}
                                                                        1 \\
                                                                        z\\
                                                                      \end{array}
                                                                    \right)-g_1(x)\left(
                                                                                    \begin{array}{c}
                                                                                      1 \\
                                                                                      0 \\
                                                                                    \end{array}
                                                                                  \right)\right),\;x\in[0,1),\nn
\ee
for every $g\in \widetilde{\mathcal{H}}$, where the Green's function $\mathcal{G}$ is given by
\be
\mathcal{G}(x,t)=\left(
                   \begin{array}{c}
                     1 \\
                     z \\
                   \end{array}
                 \right)\frac{1}{W(\psi_1,\psi_2)}\left\{\begin{array}{cc}
                                                       \psi_1(x)\psi_2(t), & t\in[0,x), \\
                                                       \psi_1(t)\psi_2(x) & t\in[x,1),
                                                     \end{array}\right.\nn
\ee
and $\psi_1,\psi_2$ are linearly independent solutions of \eqref{GIS} such that $\psi_2(0)=0$, $\psi_1\in H'[0,1)$ and $z\psi_1\in L^2([0,1);\nu)$.
\end{lem}

At the end of this subsection, we present some properties of linear relations \cite{Kato84}. For two closed linear relations $S,\ T$ in a Hilbert space $X$, the \emph{gap} is defined by
    $$\hat{\delta}(S,T):=\max\{\delta(S,T),\,\delta(T,S)\},\vspace{-0.2cm}$$
where
    $$\delta(S,T):=\sup\{{\rm dist}(\varphi,\,T):\,\varphi\in S,\,\|\varphi\|=1\}.$$
When either $S$ or $T$ is $\{0\}$, the gap is understood as $0$. Obviously, $$\hat{\delta}(S^{-1},T^{-1})=\hat{\delta}(S,T).$$

Let $T,T_n \in LR(X)$ be closed. We say that $T_n$ converges to $T$ {\it in the generalized sense} if $\hat{\delta}(T_n,T)\to0$ as $n\to\infty$.  Furthermore, if $\rho(T)\cap\rho(T_n)\neq\emptyset$ for sufficiently large $n$, we say that $T_n$ converges to $T$ {\it in the sense of norm resolvent convergence} if there exists a $z\in\rho(T)\cap\rho(T_n)$ for sufficiently large $n$ such that $\|(T_n-z)^{-1}-(T-z)^{-1}\|\to0$ as $n\to\infty$.

The following properties about the gap of linear relations were introduced in \cite{XS17}. For the concept of the norm of a linear relation, see \cite{Cross98}.

\bb{lem}\lb{l4}
Let $S,T\in LR(X)$ be closed and $A\in LR(X)$ be bounded. If $D(T)\cup D(S)\subset D(A)$
 and $A(0)\subset T(0)\cap S(0)$, then
 \be
 \hat{\delta}(S+A,T+A)\leq2(1+\|A\|^2)\hat{\delta}(S,T).\nn
 \ee
\end{lem}

\bb{lem}\lb{l5}
Let $T\in LR(X)$ be closed and $A\in LR(X)$ be  bounded. If $D(T)\subset D(A)$
 and $A(0)\subset T(0)$, then
 \be
 \hat{\delta}(T+A,T)\leq\|A\|.\nn
 \ee
\end{lem}

\bb{lem}\lb{l6}
Let $T,T_n\in LR(X)$ be closed and $\rho(T)\cap\rho(T_n)\neq\emptyset$ for sufficiently large $n$. If $T_n$ converges to $T$  in the sense of norm resolvent convergence, then $T_n$ converges to $T$ in the generalized sense.
\end{lem}

\Proof
Let $z\in\rho(T)\cap\rho(T_n)$ for sufficiently large $n$. It follows from Lemma \ref{l5} that when $n\to\infty$, we have
\beaa
&&\hat{\delta}((T_n-z)^{-1},(T-z)^{-1})\\
&=&\hat{\delta}((T_n-z)^{-1}-(T-z)^{-1}+(T-z)^{-1},(T-z)^{-1})\\
&\leq&\|(T_n-z)^{-1}-(T-z)^{-1}\|\to0.
\eeaa
Therefore, by Lemma \ref{l4}, we obtain
\be
\hat{\delta}(T_n,T)=\hat{\delta}(T_n-z+z,T-z+z)\leq2(1+|z|^2)\hat{\delta}(T_n-z,T-z)\to0. \nn
\ee
This implies that $T_n$ converges to $T$ in the generalized sense. The proof is complete.
\qed

\bb{lem}\lb{l7}{\rm\cite[Theorem 4.1]{XS17}}
Let $T\in LR(X)$ be  closed  and $\Gamma$ be a compact subset of the resolvent set $\rho(T)$. Then there exists a $\delta>0$ such that $\Gamma\subset\rho(S)$ for any closed relation $S$ with $\hat{\delta}(S,T)<\delta$.
\end{lem}
\medskip
%
%
%
%
%

\section{Main results} \lb{MaRe}

In this section, we first present some properties of solutions of \eqref{GIS}, including existence and uniqueness, as well as the bounds for these solutions, and their continuous dependence on distributions and measures. Next, we investigate the asymptotic stability  of the spectra of generalized indefinite strings, and show that the spectrum of a GIS $(1,\om,\nu)$ is upper semi-continuous with respect to $\om$ and $\nu$ in a certain sense.

%
%
%
%
%

\subsection{Properties of solutions of the differential equation (2.2)}

For $n\geq1$, let $(1,\om,\nu),(1,\om_n,\nu_n)$ be GISs where ${\mathbf w}$ and ${\mathbf w_n}$ are the anti-derivatives of $\om$ and $\om_n$. For any $d_1,d_2\in{\mathbb C}$, the basic existence and uniqueness results for \eqref{GIS} with the initial conditions
\be\lb{cond}
f(0)=d_1\;\;\; \mbox{and} \;\;\; f'(0-)=d_2,
\ee
are available; see \cite[Lemma 3.2]{EckhardtK16}. The following result establishes the existence and uniqueness of solutions to initial value problems when the initial condition is specified at the endpoint $1$.

\begin{lem}\lb{lc4}
Assume that ${\mathbf w}\in L^2[0,1)$ and $\nu$ is finite. Then for any $a\in{\mathbb C}$, there exists a unique solution $f$ of \eqref{GIS} with the initial conditions
\be\lb{cond1}
f(1-)=0\;\;\; and \;\;\; f'(1-)=a.
\ee
\end{lem}

\Proof
Since ${\mathbf w}\in L^2[0,1)$ and $\nu$ is finite, it is easy to obtain that
\be\lb{e3}
\int_0^1|{\mathbf u}|^2dx<\infty,
\ee
 where $\mathbf{u}$  is the normalized anti-derivatives of  $z\om +z^2\nu$. Consequently,
\be\lb{e4}
\int_0^1|{\mathbf u}|dx\leq\left(\int_0^1|{\mathbf u}|^2dx\right)^{\frac{1}{2}}<\infty.
\ee
Let $F=\left(
               \begin{array}{c}
                 f \\
                 f'+{\mathbf u}f\\
               \end{array}
\right)$ and $R=\left(
     \begin{array}{cc}
       -{\mathbf u} & 1 \\
       -{\mathbf u}^2 & {\mathbf u} \\
     \end{array}
   \right)$. Then the differential equation \eqref{GIS} with \eqref{cond1} is equivalent to the initial value problem
\be\lb{ec11}
F'=RF,\;\;F(1-)=\left(
                         \begin{array}{c}
                           0 \\
                           a\\
                         \end{array}
                       \right).
\ee
By \eqref{e3} and \eqref{e4}, we have
   \be\lb{e2}
   \int_0^1|R|dx=1+2\int_0^1|{\mathbf u}|dx+\int_0^1|{\mathbf u}|^2dx<\infty,
   \ee
where $|R|$ is a norm of the matrix $R$. Then  there exists a unique solution $F$ of \eqref{ec11} by \cite[Theorem 1.5.3]{Zettl}, which implies that there exists a unique solution $f$ of \eqref{GIS}--\eqref{cond1}. The proof is complete.
\qed

Now, we investigate the bounds for solutions of \eqref{GIS}.

\begin{lem}\lb{l2}
For  $d_1,d_2,a\in{\mathbb C}$, let $f$ is the solution of \eqref{GIS} with \eqref{cond} or \eqref{cond1}. Assume that ${\mathbf w}\in L^2[0,1)$ and $\nu$ is finite.   Then $f(x)$ is uniformly bounded on $[0,1]$ and
\be\lb{e1}
\int_0^{1}|f'(x)|^2dx<\infty.
\ee
\end{lem}

\Proof
Let ${\mathbf u}$, $F$, and $R$ be as in the proof of Lemma \ref{lc4}. By \eqref{e2} and \cite[Theorem 1.5.1 and Theorem 1.5.2]{Zettl}, we have that $F$ is uniformly bounded on $[0,1]$. This implies that both $f$ and $f'+{\mathbf u}f$ are uniformly bounded on $[0,1]$. Consequently, there exists a constant $M$ such that $|f(x)|<M$ for all $x\in [0.1]$ and $\int_0^1|f'+{\mathbf u}f|^2dx<\infty$. Then we have
\beaa
\int_0^1|f'(x)|^2dx&=&\int_0^1|f'+{\mathbf u}f-{\mathbf u}f|^2dx\\
&\leq&2\int_0^1|f'+{\mathbf u}f|^2dx+2\int_0^1|{\mathbf u}|^2|f|^2dx\\
&\leq&2\int_0^1|f'+{\mathbf u}f|^2dx+2M^2\int_0^1|{\mathbf u}|^2dx<\infty.
\eeaa
This concludes the proof.
\qed

Next, we show that the solutions of \eqref{GIS} with \eqref{cond} or \eqref{cond1} are continuously dependent on the distribution $\om$ and the Borel measure $\nu$ under certain conditions.

\bb{thm}\lb{sol-den}
For $d_1,d_2 \in{\mathbb C}$, let $f$ be the solution of \eqref{GIS} with \eqref{cond}. Assume that Hypotheses {\rm(1)} and {\rm(2)} hold. Then the solutions $f_n$ of the differential equation
\be\lb{e6}
 -f_n''=z\om_n f_n+z^2\nu_n f_n
\ee
with the initial conditions \eqref{cond}
converges to $f$ uniformly on $[0,1]$, and
\be\lb{e7}
\int_0^{1}|f'_n(x)-f'(x)|^2dx\to0\;\;\mbox{as}\;\;n\to\infty.
\ee
\end{thm}

\Proof Let ${\mathbf q}$ and ${\mathbf q_n}$ be the anti-derivatives of $\nu$ and $\nu_n$. By Lemma \ref{lc3}, we have
\beaa
|{\mathbf q_n}(x)-{\mathbf q}(x)|&=&\left|\int_{[0,x)}d\nu_n-\int_{[0,x)}d\nu\right|\\
&=&\left|\int_0^x(p_n-p)dt\right|\leq\int_0^1|p_n-p|dt\to0,\;\;\mbox{as}\;\;n\to\infty.
\eeaa
This implies that
\be\lb{e9}
\int_0^1|{\mathbf u_n}-{\mathbf u}|^2dx\to0,\;\;\mbox{as}\;\;n\to\infty,
\ee
and \eqref{e3} holds, where ${\mathbf u}$ and ${\mathbf u_n}$ are the anti-derivatives of $z\om+z^2\nu$ and $z\om_n+z^2\nu_n$. Thus
\be\lb{e10}
\int_0^1|{\mathbf u_n}|^2dx \leq 2\int_0^1|{\mathbf u_n}-{\mathbf u}|^2dx+2\int_0^1|{\mathbf u}|^2dx<\infty.
\ee
Let $F,R$ be as in the proof of Lemma \ref{lc4} and
$$F_n=\left(
               \begin{array}{c}
                 f_n \\
                 f'_n+{\mathbf u_n}f_n\\
               \end{array}
\right),\;\;\;R_n=\left(
     \begin{array}{cc}
       -{\mathbf u_n} & 1 \\
       -{\mathbf u_n}^2 & {\mathbf u_n} \\
     \end{array}
   \right).$$
Then  $F_n$ is a solution to the initial value problem
\be
F_n'=R_n F_n,\;\;\;F_n(0)=\left(
                         \begin{array}{c}
                           d_1 \\
                           d_2\\
                         \end{array}
                       \right).\nn
\ee
By \eqref{e3}, \eqref{e9} and \eqref{e10},  we have
   \beaa
   &&\int_0^1|R_n-R|dx\\
   &=&2\int_0^1|{\mathbf u_n}-{\mathbf u}|dx+\int_0^1|{\mathbf u_n}^2-{\mathbf u}^2|dx\\
   &\leq&2\left(\int_0^1|{\mathbf u_n}-{\mathbf u}|^2dx\right)^{\frac{1}{2}}+\left(\int_0^1|{\mathbf u_n}-{\mathbf u}|^2dx\right)^{\frac{1}{2}}\left(\int_0^1|{\mathbf u_n}|^2dx\right)^{\frac{1}{2}} \\
   &+&\left(\int_0^1|{\mathbf u_n}-{\mathbf u}|^2dx\right)^{\frac{1}{2}}\left(\int_0^1|{\mathbf u}|^2dx\right)^{\frac{1}{2}}\to0, \;\;\mbox{as}\;\;n\to\infty.
   \eeaa
It follows from \cite[Theorem 1.6.3]{Zettl} that $F_n$ converges to $F$ uniformly on $[0,1]$. Both $f_n$ and $f'_n+{\mathbf u_n}f_n$ converge to $f$ and $f'+{\mathbf u}f$ uniformly on $[0,1]$. Then
\be\lb{e11}
\int_0^1|f'_n+{\mathbf u_n}f_n-f'-{\mathbf u}f|^2dx\to0, \;\;\mbox{as}\;\;n\to\infty.
\ee
It follows from Lemma \ref{l2}, \eqref{e9} and \eqref{e10} that for any $\epsilon>0$, there exist $N,M>0$ such that when $n>N$, we have
\be
\int_0^1|{\mathbf u_n}-{\mathbf u}|^2dx<\epsilon,\;\;\;\; |f_n(x)-f(x)|^2<\epsilon,\;\;x\in[0,1),\nn
\ee
and
\be
\int_0^1|{\mathbf u_n}|^2dx<M,\;\;\;\; |f(x)|^2<M, \;\;x\in[0,1).\nn
\ee
Hence, we obtain
\beaa
&&\int_0^1|{\mathbf u_n}f_n-{\mathbf u}f|^2dx\\
&\leq& 2\int_0^1|{\mathbf u_n}f_n-{\mathbf u_n}f|^2dx+2\int_0^1|{\mathbf u_n}f-{\mathbf u}f|^2dx\\
&=&2\int_0^1|{\mathbf u_n}|^2|f_n-f|^2dx+2\int_0^1|{\mathbf u_n}-{\mathbf u}|^2|f|^2dx<4M\epsilon,\\
\eeaa
which yields that
\be
\int_0^1|{\mathbf u_n}f_n-{\mathbf u}f|^2dx\to0,\;\;\mbox{as}\;\;n\to\infty.\nn
\ee
This together with \eqref{e11} implies that
\beaa
&&\int_0^{1}|f'_n(x)-f'(x)|^2dx \\
&=&\int_0^{1}|f'_n(x)+{\mathbf u_n}f_n-{\mathbf u_n}f_n+{\mathbf u}f-{\mathbf u}f-f'(x)|^2dx\\
&\leq&2\int_0^1|f'_n+{\mathbf u_n}f_n-f'-{\mathbf u}f|^2dx+2\int_0^1|{\mathbf u_n}f_n-{\mathbf u}f|^2dx\to0.
\eeaa
Then \eqref{e7} holds. This concludes the proof.
\qed

\begin{rem} \lb{r-sol-den} {\rm By a similar proof as above, Theorem \ref{sol-den} also holds for the initial conditions \eqref{cond1}.
}
\end{rem}    
\medskip
%
%
%
%
%

\subsection{Asymptotic stability  of the spectra of generalized indefinite strings}

For $n\geq1$, let $(1,\om,\nu)$, $(1,\om_n,\nu_n)$ be GISs with the corresponding self-adjoint linear relations $T$ and $T_n$, and $\mathcal{H}_n:=H'_{*}[0,1)\times L^2([0,1);\nu_n)$. Define the operator $Q_n:\mathcal{H}_n\to \widetilde{\mathcal{H}}$ by
\be
Q_nf:=(f_1,P_nf_2),\;\;f=(f_1,f_2)\in\mathcal{H}_n,
\ee
where $P_nf_2=\sqrt{p_n}f_2$. Set $\widetilde{T}_n:=Q_nT_nQ^{-1}_n=Q_nT_nQ_n^*$. Then $\widetilde{T}_n$ is unitarily equivalent to $T_n$.

\bb{thm}\lb{main}
Assume that Hypotheses {\rm(1)} and {\rm(2)} hold. Then $\widetilde{T}_n$ converges to $\widetilde{T}$  in the sense of norm resolvent convergence.
\end{thm}

\Proof
Let $\psi_i,\,i=1,2$ be linearly independent solutions of \eqref{GIS} such that
\be\lb{3.1}
\psi_1(0)=1,\,\psi'_1(0-)=d,\,\psi_2(0)=0, \,\psi'_2(0-)=1,
\ee
with $\psi_1\in H'[0,1)$ and $z\psi_1\in L^2([0,1);\nu)$. Then $\psi_1(1-)=0$. Assume that  $\psi'_{1}(1-)=a$.
Let $\phi_n,\psi_{n,2}$ be  solutions of \eqref{e6} such that
\be\lb{3.2}
\phi_{n}(1-)=0,\,\phi'_{n}(1-)=a,\,\psi_{n,2}(0)=0, \,\psi'_{n,2}(0-)=1.
\ee
It follows from Remark \ref{r-sol-den} that $\phi_n(x)$ converges to $\psi_1(x)$ uniformly on $[0,1]$ as $n\to\infty$. Let $\psi_{n,1}:=\frac{\phi_n}{\phi_n(0)}$. Then $\psi_{n,1}(x)$ converges to $\psi_1(x)$ uniformly on $[0,1]$, and $\psi_{n,1}(1-)=0$, $\psi_{n,1}(0)=1$. Consequently,
\be
W(\psi_1,\psi_2)=W(\psi_{n,1},\psi_{n,2})=1,\nn
\ee
which implies that $\psi_{n,1}$ and $\psi_{n,2}$ are linearly independent. It follows from Theorem \ref{sol-den}, Lemma \ref{lc3} and Lemma \ref{l2}, and the integrability of $g$ that for any $\epsilon>0$, there exist $N,M>0$ such that when $n>N$, we have
\bea\lb{3.3}
&&|\psi_{n,i}(x)-\psi_i(x)|^2<\epsilon, \,\,|\psi_i(x)|^2<M,\,\,|\psi_{n,i}(x)|^2<M,\;\;x\in[0,1),\nn\\
&&\int_0^1|\psi'_{n,i}(x)-\psi'_i(x)|^2dx<\epsilon,\,\int_0^1|\psi'_{i}(x)|^2dx<M,\,\int_0^1|\psi'_{n,i}(x)|^2dx<M,\nn\\
&& \int_0^1\left|\sqrt{p_n}-\sqrt{p}\right|^2dx<\epsilon,\, \int_0^1g(x)dx<M,
\eea
and $\psi_{n,1}\in H'[0,1)$ and $z\psi_{n,1}\in L^2([0,1);\nu)$.

Let $z\in{\mathbb C}\backslash{\mathbb R}$. Then $z\in\rho(\widetilde{T})\cap\rho(\widetilde{T}_n)$ since $\widetilde{T}$ and $\widetilde{T}_n$ are self-adjoint. Set
\be
f(x,t)=\left\{\begin{array}{cc}
                                                       \psi_1(x)\psi_2(t), & t\in[0,x), \\
                                                       \psi_1(t)\psi_2(x), & t\in[x,1),
                                                       \end{array}\right.
 \; f_n(x,t)=\left\{\begin{array}{cc}
                                                       \psi_{n,1}(x)\psi_{n,2}(t), & t\in[0,x), \\
                                                       \psi_{n,1}(t)\psi_{n,2}(x), & t\in[x,1),
                                                     \end{array}\right.\nn
 \ee
and

\be
\mathcal{G}(x,t)=\left(
                   \begin{array}{c}
                     f(x,t) \\
                     zf(x,t) \\
                   \end{array}
                 \right),\;\;\;\mathcal{G}_n(x,t)=\left(
                   \begin{array}{c}
                     f_n(x,t) \\
                     zf_n(x,t) \\
                   \end{array}
                 \right).
\ee
It is easy to show that $\mathcal{G}(x,\cdot),\mathcal{G}_n(x,\cdot)\in\mathcal{H}\cap\mathcal{H}_n$ for $x\in[0,1)$.
Let $h(x)=\langle g, Q\mathcal{G}(x,\cdot)^*\rangle_{\widetilde{\mathcal H}}$ and $h_n(x)=\langle g, Q_n\mathcal{G}_n(x,\cdot)^*\rangle_{\widetilde{\mathcal H}}$ for $g\in{\widetilde{\mathcal H}}$ and $x\in[0,1)$. Then when $n>N$, by \eqref{3.3}, we have
\bea\lb{hn}
&&|h_n(x)|^2\leq\|g\|^2_{\widetilde{\mathcal H}}\|Q_n\mathcal{G}_n(x,\cdot)^*\|^2_{\widetilde{\mathcal H}}\nn\\
&=&\|g\|^2_{\widetilde{\mathcal H}}\left(|\psi_{n,1}(x)|^2\int_0^x|\psi_{n,2}'(t)|^2dt+|\psi_{n,2}(x)|^2\int_x^1|\psi_{n,1}'(t)|^2dt\right.\nn\\
&&\left.+|z|^2|\psi_{n,1}(x)|^2\int_0^xp_n|\psi_{n,2}(t)|^2dt+|z|^2|\psi_{n,2}(x)|^2\int_x^1p_n|\psi_{n,1}(t)|^2dt\right)\nn\\
&\leq&\|g\|^2_{\widetilde{\mathcal H}}\left(2M^2+2|z|^2M^2\int_0^1gdt\right)\leq\|g\|^2_{\widetilde{\mathcal H}}(2M^2+2|z|^2M^3),
\eea
which implies that $h_n\in L^2([0,1);\nu)\cap L^2([0,1);\nu_n)$. Similarly, one can show that $h\in L^2([0,1);\nu)\cap L^2([0,1);\nu_n)$. When $n>N$, by Lemma \ref{l3}, we obtain
\be\lb{3.4}
z(\widetilde{T}_n-z)^{-1}g(x)-z(\widetilde{T}-z)^{-1}g(x)=\left(
                                    \begin{array}{c}
                                      h_n(x)-h(x) \\
                                      z(P_n h_n(x)-Ph(x)) \\
                                    \end{array}
                                  \right).
\ee
Then by \eqref{3.3} and \eqref{hn}, we have
\beaa
&&\|z(\widetilde{T}_n-z)^{-1}g(x)-z(\widetilde{T}-z)^{-1}g(x)\|^2_{\widetilde{\mathcal H}}\\
&=&\int_0^1|h'_n-h'|^2dx+|z|^2\int_0^1|P_nh_n-Ph|^2dx\nn\\
&\leq&\int_0^1|h'_n-h'|^2dx+2|z|^2\int_0^1|P(h_n-h)|^2dx+2|z|^2\int_0^1|\sqrt{p_n}-\sqrt{p}|^2|h_n|^2dx\\
&\leq&\int_0^1|h'_n-h'|^2dx+2|z|^2\int_{[0,1)}|(h_n-h)|^2d\nu(x)+\|g\|^2_{\widetilde{\mathcal H}}(1+|z|^2M)4|z|^2M^2\epsilon.\\
\eeaa
Note that
\be
|h_n(x)-h(x)|^2\leq\|g\|^2_{\widetilde{\mathcal H}}\|Q_n\mathcal{G}_n(x,\cdot)^*-Q\mathcal{G}(x,\cdot)^*\|^2_{\widetilde{\mathcal H}},\nn
\ee
and
\be
|h'_n(x)-h'(x)|^2\leq\|g\|^2_{\widetilde{\mathcal H}}\left\|Q_n\frac{\partial\mathcal{G}_n}{\partial x}(x,\cdot)^*-Q\frac{\partial\mathcal{G}}{\partial x}(x,\cdot)^*\right\|^2_{\widetilde{\mathcal H}}.\nn
\ee
Then
\bea\lb{3.5}
&&\|z(\widetilde{T}_n-z)^{-1}g(x)-z(\widetilde{T}-z)^{-1}g(x)\|^2_{\widetilde{\mathcal H}}\nn\\
&\leq&\|g\|^2_{\mathcal H}\left(2|z|^2\int_{[0,1)}\|Q_n\mathcal{G}_n(x,\cdot)^*-Q\mathcal{G}(x,\cdot)^*\|^2_{\widetilde{\mathcal H}}d\nu(x)\right.\nn\\
&+&\left.\int_0^1\left\|Q_n\frac{\partial\mathcal{G}_n}{\partial x}(x,\cdot)^*-Q\frac{\partial\mathcal{G}}{\partial x}(x,\cdot)^*\right\|^2_{\widetilde{\mathcal H}}dx+(1+|z|^2M)4|z|^2M^2\epsilon\right).\nn\\
\eea
On the other hand, we have
\beaa
&&\|Q_n\mathcal{G}_n(x,\cdot)^*-Q\mathcal{G}(x,\cdot)^*\|^2_{\widetilde{\mathcal H}}=\int_0^1\left|\frac{\partial f_{n}}{\partial t}-\frac{\partial f}{\partial t}\right|^2dt\nn\\
&&+|z|^2\int_0^1|P_nf_{n}(x,t)-Pf(x,t)|^2dt:=I_1+I_2.
\eeaa
It follows from \eqref{3.3} that for $n>N$, we obtain
\beaa
I_1&=&\int_0^x|\psi_{n,1}(x)\psi'_{n,2}(t)-\psi_1(x)\psi'_2(t)|^2dt+\int_x^1|\psi'_{n,1}(t)\psi_{n,2}(x)-\psi'_1(t)\psi_2(x)|^2dt\\
&\leq&2|\psi_{n,1}(x)-\psi_1(x)|^2\int_0^x|\psi'_{n,2}(t)|^2dt+2|\psi_1(x)|^2\int_0^x|\psi'_{n,2}(t)-\psi'_{2}(t)|^2dt\\
&+&2|\psi_{n,2}(x)-\psi_2(x)|^2\int_x^1|\psi'_{n,1}(t)|^2dt+2|\psi_2(x)|^2\int_x^1|\psi'_{n,1}(t)-\psi'_{1}(t)|^2dt\\
&<&8M\epsilon.
\eeaa
Similarly, we have
\beaa
I_2&\leq&2|z|^2\int_0^1|\sqrt{p_n}-\sqrt{p}|^2|f_n|^2dt+2|z|^2\int_{[0,1)}|f_n-f|^2d\nu(t)\\
&<&4|z|^2M(M+4K)\epsilon,
\eeaa
where $K:=\nu([0,1))$. This implies that
\be\lb{3.6}
\int_{[0,1)}\|Q_n\mathcal{G}_n(x,\cdot)^*-Q\mathcal{G}(x,\cdot)^*\|^2_{\widetilde{\mathcal H}}d\nu(x)<4MK(2+|z|^2M+4|z|^2K)\epsilon.
\ee
Note that
\beaa
&&\left\|Q_n\frac{\partial\mathcal{G}_n}{\partial x}(x,\cdot)^*-Q\frac{\partial\mathcal{G}}{\partial x}(x,\cdot)^*\right\|^2_{\widetilde{\mathcal H}}=\int_0^1\left|\frac{\partial^2 f_{n}}{\partial x\partial t}-\frac{\partial^2 f}{\partial x\partial t}\right|^2dt\nn\\
&&+|z|^2\int_0^1\left|P_n\frac{\partial f_{n}}{\partial x}-P\frac{\partial f}{\partial x}\right|^2dt:=J_1+J_2.
\eeaa
Again by \eqref{3.3}, we obtain
\beaa
J_1
<2M(|\psi'_{n,1}(x)-\psi'_1(x)|^2+|\psi'_{n,2}(x)-\psi'_2(x)|^2)+2\epsilon(|\psi'_1(x)|^2+|\psi'_2(x)|^2),
\eeaa
and
\beaa
J_2&\leq&2|z|^2\int_0^1|\sqrt{p_n}-\sqrt{p}|^2\left|\frac{\partial f_n}{\partial x}\right|^2dt+2|z|^2\int_{[0,1)}\left|\frac{\partial f_n}{\partial x}-\frac{\partial f}{\partial x}\right|^2d\nu(t)\\
&<&4|z|^2MK(|\psi'_{n,1}(x)-\psi'_1(x)|^2+|\psi'_{n,2}(x)-\psi'_2(x)|^2)\\
&+&4|z|^2K\epsilon(|\psi'_1(x)|^2+|\psi'_2(x)|^2)+2|z|^2M\epsilon(|\psi'_{n,1}(x)|^2+|\psi'_{n,2}(x)|^2).
\eeaa
Thus
\bea\lb{3.7}
&&\int_0^1\left\|Q_n\frac{\partial\mathcal{G}_n}{\partial x}(x,\cdot)^*-Q\frac{\partial\mathcal{G}}{\partial x}(x,\cdot)^*\right\|^2_{\widetilde{\mathcal H}}dx\nn\\
&<&2M(1+2|z|^2K)\int_0^1(|\psi'_{n,1}(x)-\psi'_1(x)|^2+|\psi'_{n,2}(x)-\psi'_2(x)|^2)dx\nn\\
&+&2|z|^2M\epsilon\int_0^1(|\psi'_{n,1}(x)|^2+|\psi'_{n,2}(x)|^2)dx\nn\\
&+&2\epsilon(1+2|z|^2K)\int_0^1(|\psi'_1(x)|^2+|\psi'_2(x)|^2)dx\nn\\
&<&4M(2+|z|^2M+4|z|^2K)\epsilon.
\eea
It follows from \eqref{3.5}-\eqref{3.7} that when $n>N$, we have
\beaa
&&\|z(\widetilde{T}_n-z)^{-1}-z(\widetilde{T}-z)^{-1}\|^2_{\widetilde{\mathcal H}}\\
&<&8M\epsilon+8|z|^2M(M+4K)\epsilon+4|z|^4M(M^2+2MK+8K^2)\epsilon.
\eeaa
This yields that $\|(\widetilde{T}_n-z)^{-1}-(\widetilde{T}-z)^{-1}\|_{\widetilde{\mathcal H}}\to0$ as $n\to\infty$, that is, $\widetilde{T}_n$ converges to $\widetilde{T}$  in the sense of norm resolvent convergence. The proof is complete.
\qed

\bb{rem}\lb{c1}
{\rm By Lemma \ref{l6} and Theorem \ref{main}, we obtain that $\widetilde{T}_n$ converges to $\widetilde{T}$  in the generalized sense.}
\end{rem}

Since $T$ and $T_n$ are unitarily equivalent to $\widetilde{T}$ and $\widetilde{T}_n$, the following result  is deduced from Lemma \ref{l7} and Remark \ref{c1}.

\bb{cor}\lb{reso} Assume that Hypotheses {\rm(1)} and {\rm(2)} hold. If $\Gamma$ is a compact subset of the resolvent set $\rho(T)$, then  $\Gamma\subset \rho(T_n)$ for sufficiently large $n$.
\end{cor}

The following result shows that the spectrum of a GIS is upper semi-continnous on $\om$ and $\nu$  under Hypotheses {\rm(1)} and {\rm(2)}. We use the notation in \cite[p. 208]{Kato84} as follows
\be
\varrho(\sigma(T_n),\sigma(T)):=\sup_{\lambda\in\sigma(T_n)}{\rm dist}(\lambda,\sigma(T)).\nn
\ee

\bb{thm}\lb{tt}  Assume that Hypotheses {\rm(1)} and {\rm(2)} hold. Then $\varrho(\sigma(T_n),\sigma(T))\to0$ as $n\to\infty$.
\end{thm}

\Proof
For any $\epsilon>0$, let
\be
\Gamma:=\{\lambda\in{\mathbb C}:{\rm dist}(\lambda,\sigma(T))\geq\varepsilon\}.\nn
\ee
Then $\Gamma$ is a compact subset of $\rho(T)$. It follows from Corollary \ref{reso} that there exist a $N>0$ such that when $n>N$, we have $\Gamma\subset \rho(T_n)$. It follows that
\be
\sigma(T_n)\subset\{\lambda\in{\mathbb C}:{\rm dist}(\lambda,\sigma(T))<\varepsilon\}.\nn
\ee
Thus we have $\varrho(\sigma(T_n),\sigma(T))\leq\epsilon$, which implies that $\varrho(\sigma(T_n),\sigma(T))\to 0$ as $n\to \infty$. This concludes the proof.
\qed

\bb{rem}
{\rm The convergence of distributions in Definition \ref{dis-uc} is stronger than the weak$^*$ convergence of distributions; see \cite{Grubb09}. A natural question is to study the behaviour of the spectrum of $T_n$ if  $\omega_n$ weak$^*$ converges to $\omega$. We conjecture that  Theorem \ref{tt} is true if the Hypothesis {\rm(2)} is replaced by the condition that $\omega_n$ weak$^*$ converges to $\omega$, and we plan to discuss this issue in further work as did in \cite{Chu20, Meng13}.
}
\end{rem}

\bb{rem}
{\rm The differential equation \eqref{GIS} is equivalent to the differential equation
 \be\lb{GIS'}
 -f''=z\om f+z^2pf,
 \ee
 where $\nu$ is a positive finite Borel measures on $[0,1)$ that is absolutely continuous with respect to the Lebesgue measure, and $p\in L^1[0,1)$ is the corresponding Radon-Nikodym derivative. Therefore, the results in this paper also can be applied to \eqref{GIS'}.
 }
\end{rem}
\medskip
%
%
%
%
%
\section*{Acknowledgments}
The authors would like to thank Professor Jifeng Chu for helpful suggestions.

\section*{Date Availability}
No data were used to support this study.

\section*{Conflicts of Interest}
The authors have no competing interests to declare that are relevant to the content of this article.

\end{document}